\documentclass{article}

\usepackage[english]{babel}

\usepackage[letterpaper,top=2cm,bottom=2cm,left=3cm,right=3cm,marginparwidth=1.75cm]{geometry}

\usepackage{amsmath}
\usepackage{amssymb}
\usepackage{amsthm}
\usepackage{graphicx}
\usepackage[colorlinks=true, allcolors=blue]{hyperref}

\theoremstyle{plain}
\newtheorem{theorem}{Theorem}
\newtheorem{lemma}[theorem]{Lemma}
\newtheorem{proposition}[theorem]{Proposition}
\newtheorem{conjecture}[theorem]{Conjecture}
\newtheorem{corollary}[theorem]{Corollary}
\newtheorem{problem}[theorem]{Problem}

\theoremstyle{definition}
\newtheorem{definition}[theorem]{Definition}
\newtheorem{example}[theorem]{Example}

\title{Asymmetric SICs over finite fields}
\author{Joseph W.\ Iverson\thanks{Department of Mathematics, Iowa State University, Ames, IA} \and Dustin G.\ Mixon\thanks{Department of Mathematics, The Ohio State University, Columbus, OH} \thanks{Translational Data Analytics Institute, The Ohio State University, Columbus, OH}}
\date{}

\begin{document}
\maketitle

\begin{abstract}
Zauner's conjecture concerns the existence of $d^2$ equiangular lines in $\mathbb{C}^d$; such a system of lines is known as a SIC.
In this paper, we construct infinitely many new SICs over finite fields.
While all previously known SICs exhibit Weyl--Heisenberg symmetry, some of our new SICs exhibit trivial automorphism groups.
We conjecture that such \textit{totally asymmetric} SICs exist in infinitely many dimensions in the finite field setting.
\end{abstract}

\section{Introduction}

Given a field $K$ over which the polynomial $t^2+1$ does not split (so that, in particular, $\operatorname{char} K \neq 2$), consider any splitting field extension $L$ obtained by adjoining a root $\mathrm{i}$ of this polynomial.
For example, if $K=\mathbb{R}$, then $L=\mathbb{C}$, while if $K=\mathbb{F}_p$ for some prime $p\equiv 3\bmod 4$, then $L\cong\mathbb{F}_{p^2}$. 
The \textit{conjugate} of $a+b\mathrm{i}\in L$ with $a,b \in K$ is then defined as $\overline{a+b\mathrm{i}}=a-b\mathrm{i}$, and given a column vector $x\in L^d$, we let $x^*$ denote its conjugate transpose.
This determines a sesquilinear form $(\cdot,\cdot)$ on $L^d$ defined by $(x,y)=x^*y$.
Throughout, we use the notation $[d]:=\{1,\ldots,d\}$, we let $I$ denote the $d\times d$ identity matrix, we take $e_i$ to be the $i$th column of $I$, and we denote the all-ones vector by $\mathbf{1}=\sum_{i\in[d]}e_i$.
What follows is our primary object of interest.

\begin{definition}
\label{defn:SIC}
Given $a,b,c\in K$ such that $a^2\neq b$,
and given a set $X$ of cardinality $d^2$, the tuple $\{ x_i \}_{i\in X}$ of vectors $x_i \in L^d$ is said to form an \textbf{$(a,b,c)$-SIC} if
\begin{itemize}
\item[(a)]
$(x_i,x_i)=a$ for all $i\in X$,
\item[(b)]
$(x_i,x_j)(x_j,x_i)=b$ for all $i,j\in X$ with $i\neq j$,
\item[(c)]
$\sum_{i\in X}x_ix_i^* = cI$, and
\item[(d)]
$\operatorname{span}\{ x_i : i \in X \} = L^d$.
\end{itemize}
\end{definition}

The name ``SIC'' originates from the case $L=\mathbb{C}$, where SIC is short for SIC-POVM, which in turn is short for the quantum information--theoretic jargon \textit{symmetric, informationally complete positive operator--valued measure}~\cite{RenesBKSC:04}.
Note that in this complex case, SIC vectors span points in projective space that are pairwise equidistant in the Fubini--Study metric.
In fact, complex SICs achieve equality in Gerzon's upper bound on the size of such configurations~\cite{LemmensS:73}.
Zauner conjectured that complex SICs exist in every dimension~\cite{Zauner:99}.
To date, there are only finitely many dimensions for which complex SICs are known to exist, though Zauner's conjecture was recently shown to be a consequence of conjectures in algebraic number theory~\cite{ApplebyFK:25}.
Interestingly, Gerzon's bound also holds over finite fields (see Theorem~4.2 in~\cite{GreavesIJM:22}), where infinitely many SICs have been discovered~\cite{GreavesIJM:22,IversonKM:21}.

To date, every known SIC (over any field $L$) can be constructed (up to trivial equivalences) by first selecting an abelian group $G$ of order $d$, and then taking all translations and modulations of a specially selected seed function $G\to L$, resulting in a system of vectors that exhibits \textit{Weyl--Heisenberg} symmetry.
In this paper, we present the first known SICs that do not exhibit such symmetry.
In the next section, we introduce our general construction and consider several examples to illustrate the variety of SICs we construct.
Next, Section~3 studies the (lack of) symmetry in these SICs.
We conclude in Section~4 with examples of asymmetric SICS and a discussion of implications for Zauner's conjecture.

\section{A fruitful construction}

What follows is our main result, which uses (modular) Hadamard matrices to generalize Hoggar's lines~\cite{Hoggar:81,Hoggar:98}, similar to~\cite{JedwabW:15}.

\begin{theorem}
\label{thm.main}
Fix a dimension $d\equiv 8\bmod\operatorname{char}K$.
Select any modular Hadamard matrix $H\in K^{d\times d}$, meaning $H_{ij}^2=1$ for all $i,j\in[d]$ and $H^\top H=dI$, and let $h_j\in K^{d}$ denote the $j$th column of $H$.
Then the vectors $\{x_{ij}\}_{i,j\in[d]}$ in $L^d$ defined in terms of the Hadamard product by
\[
x_{ij}
=h_j\circ(\mathbf{1}+ze_i),
\qquad
z:=-2(1+\mathrm{i})
\]
form a $(12,16,96)$-SIC.
\end{theorem}

This construction can be viewed as a combinatorial analog of taking all translations and modulations of a seed function.
In particular, if we interpret the index set $[d]$ as the carrier set of an abelian group $G$, then the vectors $\{\mathbf{1}+ze_i\}_{i\in[d]}$ form an orbit under the regular representation of $G$.
Next, instead of modulating each orbit representative with the characters of $G$, we use the columns of a (modular) Hadamard matrix.
A similar construction appeared in~\cite{JedwabW:15}, which instead used complex Hadamard matrices and considered all possible choices $z$, but the construction only produced complex SICs in dimensions $2$, $3$, and $8$.
By contrast, we obtain a multitude of SICs over finite fields, as elucidated further below.
To prove Theorem~\ref{thm.main}, we consider the Gram matrix of $\{x_{ij}\}_{i,j\in[d]}$, which is given by the following.

\begin{lemma}
\label{lem.gram matrix}
The vectors $\{x_{ij}\}_{i,j\in[d]}$ constructed in Theorem~\ref{thm.main} satisfy
\begin{equation}
\label{eq:gram matrix}
(x_{ij},x_{kl})
=4\cdot(-1)^{\delta_{ik}}\cdot(-1)^{\delta_{jl}}\cdot\left\{
\begin{array}{rl}
+1&\text{if } ~ H_{ij}H_{il}=-1 ~ \text{ and } ~ H_{kj}H_{kl}=-1\\
-1&\text{if } ~ H_{ij}H_{il}=+1 ~ \text{ and } ~ H_{kj}H_{kl}=+1\\
+\mathrm{i}&\text{if } ~ H_{ij}H_{il}=+1 ~ \text{ and } ~ H_{kj}H_{kl}=-1\\
-\mathrm{i}&\text{if } ~ H_{ij}H_{il}=-1 ~ \text{ and } ~ H_{kj}H_{kl}=+1
\end{array}
\right.
\end{equation}
whenever $(i,j)\neq(k,l)$.
\end{lemma}

\begin{proof}
It is convenient to write $x_{ij}=h_j+zH_{ij}e_i$.
Then
\[
(x_{ij},x_{kl})
=(h_j+zH_{ij}e_i,h_l+zH_{kl}e_k)
=d\delta_{jl}+zH_{kj}H_{kl}+\overline{z}H_{ij}H_{il}+z\overline{z}H_{ij}H_{kl}\delta_{ik}.
\]
We proceed in cases.

\medskip

\noindent
\textbf{Case I:} $i\neq k$, $j\neq l$.
Then
\begin{align*}
(x_{ij},x_{kl})
=zH_{kj}H_{kl}+\overline{z}H_{ij}H_{il}
&=-2(1+\mathrm{i})H_{kj}H_{kl}-2(1-\mathrm{i})H_{ij}H_{il}\\
&=-2(H_{kj}H_{kl}+H_{ij}H_{il})-2\mathrm{i}(H_{kj}H_{kl}-H_{ij}H_{il}),
\end{align*}
which reduces further to the right-hand side of~\eqref{eq:gram matrix}.

\medskip
\noindent
\textbf{Case II:} $i\neq k$, $j=l$.
Then the right-hand side of~\eqref{eq:gram matrix} equals~4, and indeed,
\[
(x_{ij},x_{kl})
=d+zH_{kj}H_{kj}+\overline{z}H_{ij}H_{ij}
=d+z+\overline{z}
=8-4
=4.
\]

\medskip
\noindent
\textbf{Case III:} $i=k$, $j\neq l$.
Then the right-hand side of~\eqref{eq:gram matrix} equals $4 H_{ij} H_{il}$, and indeed,
\[
(x_{ij},x_{kl})
=zH_{kj}H_{kl}+\overline{z}H_{ij}H_{il}+z\overline{z}H_{ij}H_{kl}
=H_{ij}H_{il}(z+\overline{z}+z\overline{z})
=4H_{ij}H_{il}.
\qedhere
\]
\end{proof}

\begin{proof}[Proof of Theorem~\ref{thm.main}]
We demonstrate each property from Definition~\ref{defn:SIC}.
For (a), we have
\[
(x_{ij},x_{ij})
=(d-1)+(1+z)\overline{(1+z)}
=d+z+\overline{z}+z\overline{z}
=d+4
=12.
\]
For (b), Lemma~\ref{lem.gram matrix} implies that $(x_{ij},x_{kl})(x_{kl},x_{ij})=16$ whenever $(i,j)\neq(k,l)$.
For (c), first note that since $d\not\equiv 0\bmod\operatorname{char}K$ by assumption, it holds that $d^{-1}H^\top$ is the matrix inverse of $H$, and so $HH^\top=dI$.
Next, note that
\[
x_{ij}x_{ij}^*
=(h_j+zH_{ij}e_i)(h_j+zH_{ij}e_i)^*
=h_jh_j^\top+\overline{z}H_{ij}h_je_i^\top+zH_{ij}e_ih_j^\top+8e_ie_i^\top,
\]
and so
\begin{align*}
e_k^\top(x_{ij}x_{ij}^*)e_l
&=e_k^\top(h_jh_j^\top+\overline{z}H_{ij}h_je_i^\top+zH_{ij}e_ih_j^\top+8e_ie_i^\top)e_l\\
&=H_{kj}H_{lj}+\overline{z}H_{ij}H_{kj}\delta_{il}+zH_{ij}\delta_{ik}H_{lj}+8\delta_{ik}\delta_{il}.
\end{align*}
Thus, if $k\neq l$, we have
\begin{align*}
e_k^\top\bigg(\sum_{i,j\in[d]}x_{ij}x_{ij}^*\bigg)e_l
&=\sum_{i,j\in[d]}\big(H_{kj}H_{lj}+\overline{z}H_{ij}H_{kj}\delta_{il}+zH_{ij}\delta_{ik}H_{lj}\Big)\\
&=d\sum_{j\in [d]}H_{kj}H_{lj}+\overline{z}\sum_{j\in[d]}H_{lj}H_{kj}+z\sum_{j\in[d]}H_{kj}H_{lj}
=0,
\end{align*}
while if $k=l$, we instead have
\begin{align*}
e_k^\top\bigg(\sum_{i,j\in[d]}x_{ij}x_{ij}^*\bigg)e_l
&=\sum_{i,j\in[d]}\big(1+\overline{z}H_{ij}H_{kj}\delta_{ik}+zH_{ij}\delta_{ik}H_{kj}+8\delta_{ik}\big)\\
&=d^2+d\overline{z}+dz+8d
=d^2+4d
=96,
\end{align*}
as claimed.
Finally, for (d), observe that $\left[ x_{11} \ \cdots \ x_{1d} \right] = \operatorname{diag}(1+z,1,\ldots,1) H$ already has full rank.
\end{proof}

To illustrate the significance of Theorem~\ref{thm.main}, we consider several instances of this result.
First, we interpret Theorem~\ref{thm.main} in the classical setting where $L=\mathbb{C}$.

\begin{example}
\label{ex:Hoggar}
When $L=\mathbb{C}$, Theorem~\ref{thm.main} only constructs SICs in dimension $d=8$.
In the special case where $H$ is the character table of $(\mathbb{Z}/2\mathbb{Z})^3$, this construction is known in the literature as \textit{Hoggar's lines}~\cite{Hoggar:81,Hoggar:98}.
It is known that every Hadamard matrix of order $8$ can be realized by applying appropriate signed permutations to the rows and columns of this particular Hadamard matrix, and so, as a consequence of Theorem~\ref{thm:sandwich} below, all of the SICs from Theorem~\ref{thm.main} in this case are equivalent.
(We will revisit Hoggar's lines to discuss their symmetries in Example~\ref{ex:Hoggar syms}.)
\end{example}

In the case where $L$ is a finite field, we can still use honest Hadamard matrices in Theorem~\ref{thm.main}, which are conjectured to exist whenever $d$ is a multiple of $4$.
(Also, this conjecture is currently known to hold for all $d<666$.)
What follows is an immediate consequence of Theorem~\ref{thm.main}.

\begin{corollary}
Given a Hadamard matrix of order $d$, then for every prime $p\equiv 3\bmod 4$ that divides $d-8$, there exists a $d$-dimensional SIC over $\mathbb{F}_{p^2}$.
\end{corollary}

Next, we use known parameterized families of Hadamard matrices to obtain a further consequence that was previously only known to hold when $p=3$ (see~\cite{GreavesIJM:22}).

\begin{corollary}
\label{cor.inf dim for each p}
For every prime $p \equiv 3 \bmod 4$, SICs exist in infinitely many dimensions over $\mathbb{F}_{p^2}$.
\end{corollary}

\begin{proof}
We proceed in cases.

\medskip

\noindent
\textbf{Case I:}
$p\neq 7$.
Take $d=q+1$ for any prime $q\equiv 7\bmod 4p$.
Since $q\equiv 3\bmod 4$, Paley's construction I gives a Hadamard matrix of order $d$.
Furthermore, $q\equiv 7\bmod p$, i.e., $d\equiv 8\bmod p$, and so Theorem~\ref{thm.main} delivers a $d$-dimensional SIC over $\mathbb{F}_{p^2}$.
Since $7$ and $4p$ are coprime by assumption, there are infinitely many choices for $q$ by Dirichlet's theorem on arithmetic progressions.

\medskip

\noindent
\textbf{Case II:}
$p=7$.
Take $d=2(q+1)$ for any prime $q\equiv 17\bmod 28$.
Since $q\equiv 1\bmod 4$, Paley's construction II gives a Hadamard matrix of order $d$.
Furthermore, $q+1\equiv 18\bmod 28$, meaning $d\equiv 36\bmod 56$, and so $d\equiv 8\bmod 7$.
Thus, Theorem~\ref{thm.main} delivers a $d$-dimensional SIC over $\mathbb{F}_{7^2}$.
Since $17$ and $28$ are coprime, there are infinitely many choices of $q$ by Dirichlet's theorem on arithmetic progressions.
\end{proof}

Next, instead of honest Hadamard matrices, we use $p$-modular Hadamard matrices, which in turn deliver many more SICs.

\begin{example}
When $p=3$, it is known that $p$-modular Hadamard matrices of order $d$ exist precisely when $d\not\equiv 5\bmod 6$ (see~\cite{LeeS:14}).
Considering the requirement in Theorem~\ref{thm.main} that $d\equiv 8\bmod p$, we therefore get $d$-dimensional SICs over $\mathbb{F}_{3^2}$ whenever $d\equiv 2\bmod 6$.
(Note that the proof of Corollary~\ref{cor.inf dim for each p} accounts for only a tiny portion of these dimensions.)
This directly generalizes Theorem~6.5 in~\cite{GreavesIJM:22}, which constructs SICs over $\mathbb{F}_{3^2}$ in dimension $d=2^{2k+1}$, which is also $2\bmod 6$.
\end{example}

\begin{example}
When $p=7$, it is known that $p$-modular Hadamard matrices of order $d>29$ exist when $d\equiv 1\bmod 7$ (see~\cite{Kuperberg:16}).
Since this implies $d\equiv 8\bmod p$, we therefore get $d$-dimensional SICs over $\mathbb{F}_{7^2}$ in all such cases.
Meanwhile, the proof of Corollary~\ref{cor.inf dim for each p} accounts for a tiny portion of these dimensions.
\end{example}


For larger primes (e.g., $p=11$), we turn to asymptotic existence results as in~\cite{Kuperberg:16}.
Whenever $p$ is an odd prime, a $p$-modular Hadamard matrix of order $d$ exists only if $d$ is even or a quadratic residue modulo $p$ (see~\cite{LeeS:14}), and this necessary condition is conjectured to be asymptotically sufficient:

\begin{conjecture}[The $p$-modular Hadamard conjecture; Conjecture~5.1 in~\cite{Kuperberg:16}]
Given an odd prime $p$, then for all but finitely many positive integers $d$, there exists a $p$-modular Hadamard matrix of order $d$ if and only if $d$ is even or a quadratic residue modulo $p$.
\end{conjecture}

As we will see, this conjectured plethora of modular Hadamard matrices implies that Theorem~\ref{thm.main} is particularly fruitful.
In what follows, the \textit{density} of $S\subseteq\mathbb{N}$ is defined by
\[
\operatorname{dens}(S)
=\limsup_{n\to\infty}\frac{|S\cap\{1,\ldots,n\}|}{n},
\]
and we say a set contains \textit{almost every} positive integer if its complement in $\mathbb{N}$ has density zero. 

\begin{corollary}
Conditioned on the $p$-modular Hadamard conjecture, the following holds:
For almost every positive integer $d$, there exists a SIC in $L^d$ for some finite field $L$.
\end{corollary}

\begin{proof}
For each prime $p\equiv 7\bmod 8$, let $S_p\subseteq\mathbb{N}$ denote the set of $d\equiv 8\bmod p$ for which there exists a $d\times d$ Hadamard matrix modulo $p$.
For each $d\in S_p$, we may combine this modular Hadamard matrix with Theorem~\ref{thm.main} to obtain a SIC in $L^d$ with $L=\mathbb{F}_{p^2}$.
Notice that $8$ is a quadratic residue modulo each prime $p\equiv 7\bmod 8$, as a consequence of Euler's criterion.
Thus, by the $p$-modular Hadamard conjecture, for each prime $p\equiv7\bmod8$, there exists a finite set $F_p\subseteq\mathbb{N}$ such that $S_p\cup F_p$ contains every $d\equiv 8\bmod p$.
As such, there exist $c,C>0$ such that for every sufficiently large $x$, it holds that
\begin{align*}
\operatorname{dens}\bigg(\mathbb{N}\setminus \bigcup_{\substack{p\text{ prime}\\p\equiv7\bmod8}}S_p\bigg)
&\leq\operatorname{dens}\bigg(\mathbb{N}\setminus \bigcup_{\substack{p\text{ prime}\\p\equiv7\bmod8\\p\leq x}}S_p\bigg)\\
&=\operatorname{dens}\bigg(\mathbb{N}\setminus \bigcup_{\substack{p\text{ prime}\\p\equiv7\bmod8\\p\leq x}}(S_p\cup F_p)\bigg)\\
&\leq\operatorname{dens}\bigg(\mathbb{N}\setminus \bigcup_{\substack{p\text{ prime}\\p\equiv7\bmod8\\p\leq x}}(8+p\mathbb{N})\bigg)
\leq c\prod_{\substack{p\text{ prime}\\p\equiv7\bmod8\\p\leq x}}\Big(1-\frac{1}{p}\Big)
\leq\frac{C}{(\log x)^{1/4}},
\end{align*}
where the last steps follow from a standard inclusion--exclusion argument and Mertens' third theorem for arithmetic progressions~\cite{Williams:74}.
The result then follows by taking $x\to\infty$.
\end{proof}

\section{Equivalence and symmetries}

It is generally believed that the Hadamard conjecture holds in a strong sense, namely, that the number of inequivalent Hadamard matrices of order $4t$ grows rapidly as a function of $t$.
Furthermore, it is believed that most large Hadamard matrices exhibit few symmetries. 
For example, there are already over $13$ million inequivalent Hadamard matrices of order $32$, and 
interestingly, the vast majority of these matrices have small automorphism groups~\cite{KharaghaniT:13}.
In this section, we study the relationship between the symmetries of a Hadamard matrix and the symmetries of the corresponding SIC constructed by Theorem~\ref{thm.main}.
The following definitions make precise the notions of symmetry that we consider.

\begin{definition}
Let $H$ and $H'$ be matrices in $K^{X \times X}$.
\begin{itemize}
\item[(a)]
An ordered pair of permutations $(\pi, \sigma) \in S_X \times S_X$ induces a \textbf{weak equivalence} from $H$ to $H'$ if there exist $\varepsilon_i, \varepsilon'_i \in \{\pm 1\}$ such that $H_{\pi(i),\sigma(j)}' = \varepsilon_i \varepsilon'_j H_{ij}$ for every $i,j \in X$.

\item[(b)]
A permutation $\pi \in S_X$ induces a \textbf{strong equivalence} from $H$ to $H'$ if there exist $\varepsilon_i \in \{\pm 1\}$ such that $H_{\pi(i),\pi(j)}' = \varepsilon_i \varepsilon_j H_{ij}$ for every $i,j \in X$.

\item[(c)]
The \textbf{weak automorphism group} of $H$ is the group $\operatorname{Aut}_w(H) \leq S_X \times S_X$ of all ordered pairs of permutations that induce weak equivalence from $H$ to itself.

\item[(d)]
The \textbf{strong automorphism group} of $H$ is the group $\operatorname{Aut}_s(H) \leq S_X$ of permutations that induce strong equivalence from $H$ to itself.

\end{itemize}
\end{definition}

If matrices are strongly equivalent, then they are weakly equivalent.
Furthermore, the diagonal inclusion $S_X \hookrightarrow S_X \times S_X$ embeds $\operatorname{Aut}_s(H)$ as a subgroup of $\operatorname{Aut}_w(H)$.

The weak automorphism group of $H$ is closely related to its classical \textit{automorphism group}, namely, the group of all ordered pairs $(P,Q)$ of monomial $\{0,\pm 1\}$ matrices that satisfy $PHQ^\top = H$.
To see this, and to better understand weak equivalence, consider the action $(P,Q) \cdot H = PHQ^\top$ of ordered pairs of monomial $\{0,\pm1\}$-matrices on $K^{X \times X}$.
When $(P,Q)\cdot H = H'$, the underlying pair of permutations induces a weak equivalence from $H$ to $H'$.
In particular, the weak automorphism group of $H$ consists of all pairs of permutations underlying its stabilizer.
Thus, the classical automorphism group of $H$ covers $\operatorname{Aut}_w(H)$ by sending a monomial matrix to its permutation; this is a 2-fold covering if $H$ has all nonzero entries.
Likewise, strong equivalence and the strong automorphism group can be understood in terms of the action $P \cdot H = PHP^\top$ of monomial $\{0,\pm1\}$-matrices on $K^{X \times X}$.
In particular, strongly equivalent matrices are similar matrices.

The distinction between strong and weak equivalence can be seen in the following example.

\begin{example}
There are eight $2\times 2$ Hadamard matrices, and they are all weakly equivalent with weak automorphism group $S_2 \times S_2$.
However, there are exactly three strong equivalence classes (distinguished by trace) with representatives:
\[
H_1 = \left[ \begin{array}{rr} 1 & 1 \\ 1 & -1 \end{array} \right],
\qquad
H_2 = \left[ \begin{array}{rr} 1 & 1 \\ -1 & 1 \end{array} \right],
\qquad
H_3 = \left[ \begin{array}{rr} -1 & -1 \\ 1 & -1 \end{array} \right].
\]
Here, $\operatorname{Aut}_s(H_1)$ is trivial, while $\operatorname{Aut}_s(H_2) =\operatorname{Aut}_s(H_3) = S_2$.
\end{example}

Next, we consider notions of equivalence and symmetry for SICs.

\begin{definition}
\label{defn:SIC equivalence}
Let $\{x_i\}_{i\in X}$ and $\{x_i'\}_{i\in X}$ be tuples of vectors in $L^d$.
\begin{itemize}
\item[(a)]
A permutation $\pi \in S_X$ induces a \textbf{weak equivalence} from $\{ x_i \}$ to $\{x_i'\}$ if there exist
    \begin{itemize}
    \item[(i)]
    an invertible linear map $U\colon L^d\to L^d$ that preserves the sesquilinear form $(\cdot,\cdot)$,
    \item[(ii)]
    a tuple $\{c_i\}_{i \in X}$ of scalars in $L$ for which $\overline{c_i}c_i$ is constant, and
    \item[(iii)]
    a field automorphism $\gamma \colon L \to L$ satisfying $(\overline{\alpha})^\gamma = \overline{(\alpha^\gamma)}$ for every $\alpha \in L$
    \end{itemize}
such that $x_{\pi(i)}'=c_i(Ux_i)^\gamma$ for every $i \in X$.

\item[(b)]
A permutation $\pi \in S_X$ induces a \textbf{strong equivalence} from $\{ x_i \}$ to $\{x_i'\}$ if there exist
    \begin{itemize}
    \item[(i)]
    an invertible linear map $U\colon L^d\to L^d$ that preserves the sesquilinear form $(\cdot,\cdot)$, and
    \item[(ii)]
    a tuple $\{c_i\}_{i \in X}$ of scalars in $L$ for which $\overline{c_i}c_i$ is constant
    \end{itemize}
such that $x_{\pi(i)}'=c_iUx_i$ for every $i \in X$.

\item[(c)]
The \textbf{weak automorphism group} of $\{x_i\}$ is the group $\operatorname{Aut}_w(\{x_i\}) \leq S_X$ of all permutations that induce weak equivalence from $\{x_i\}$ to itself.

\item[(d)]
The \textbf{strong automorphism group} of $\{x_i\}$ is the group $\operatorname{Aut}_s(\{x_i\}) \leq S_X$ of all permutations that induce strong equivalence from $\{x_i\}$ to itself.

\end{itemize}
\end{definition}

\begin{example}
Consider the SIC $\{x_i\}_{i=1}^4$ in $\mathbb{C}^2$ given by
\[
\left[ \begin{array}{cccc} x_1 & x_2 & x_3 & x_4 \end{array} \right]
= \left[ \begin{array}{cccc}
1 & \mu & \mu & \mu \\
0 & \sqrt{1-\mu^2} & \omega\sqrt{1-\mu^2} & \overline{\omega}\sqrt{1-\mu^2}
\end{array} \right],
\qquad
\mu = \tfrac{1}{\sqrt{3}},
\qquad
\omega = e^{2\pi i/3}.
\]
Applying complex conjugation demonstrates that $(3\ 4) \in \operatorname{Aut}_w(\{x_i\})$.
However, one can show that $(3\ 4) \notin \operatorname{Aut}_s(\{x_i\})$.
Indeed, a calculation yields
\[
(x_1, x_2)(x_2, x_3)(x_3, x_1) = \mu^3 \mathrm{i}
\qquad
\text{and}
\qquad
(x_1, x_2)(x_2, x_4)(x_4, x_1) = -\mu^3 \mathrm{i},
\]
while the argument of any such ``triple product'' is held invariant by operations of the form $x_i \mapsto c_i U x_i$ with $U \in \mathbb{C}^{2 \times 2}$ unitary and $\overline{c_i} c_i = c$ constant:
\[
(c_i Ux_i,c_j Ux_j)(c_j Ux_j,c_k Ux_k)(c_k U x_k, c_i U x_i)
= c^3(x_i,x_j)(x_j,x_k)(x_k,x_i).
\]
Hence, no such operation swaps $x_3$ and $x_4$ while fixing $x_1$ and~$x_2$.
\end{example}

In the following theorem, we relate the symmetries of the SIC in Theorem~\ref{thm.main} with those of the modular Hadamard matrix employed in its construction, as well as those of a related matrix.
Specifically, given a modular Hadamard matrix $H\in K^{d\times d}$ whose $i$th row is $r_i^\top$, we denote $\tilde{H} \in K^{[d]^2 \times [d]^2} \cong (K^{d\times d})^{d\times d}$ for the modular Hadamard matrix whose $(i,k)$th block is given by $r_k r_i^\top$; that is, $\tilde{H}_{(i,j),(k,l)} = [r_k r_i^\top]_{jl}$.
This construction is well known; for instance, it appears in~\cite{HolzmannKO:10}.

The significance of $\tilde{H}$ in this setting is as follows.
Suppose Theorem~\ref{thm.main} creates a SIC $\{x_{ij}\}_{i,j \in [d]}$ from $H$.
Then the tensors $x_{ij} \otimes x_{ij}$ span equiangular lines in $L^d \otimes L^d$; in fact, they form an \textit{equiangular tight frame} for the space of symmetric tensors when $K = \mathbb{R}$ or $K$ is finite~\cite{Waldron:17,IversonKM:21}.
Changing line representatives, one can show that $\{ H_{ij}\, x_{ij} \otimes x_{ij} \}_{i,j \in [d]}$ has Gram matrix $16\tilde{H}+128I$.
Given this correspondence, one might expect the symmetries of $\{x_{ij}\}$ to manifest as symmetries of $\tilde{H}$.
Indeed, we have the following.

\begin{theorem}[Automorphism group sandwich]
\label{thm:sandwich}
Given modular Hadamard matrices $H$ and $H'$ in $K^{d \times d}$, consider the corresponding SICs $\{x_{ij}\}_{(i,j) \in [d]^2}$ and $\{x_{ij}'\}_{(i,j) \in [d]^2}$ in $L^d$ constructed by Theorem~\ref{thm.main}.
Next, 
let $\iota \colon S_d \times S_d \to S_{[d]\times [d]}$ be the embedding given by 
\[
\iota(\pi,\sigma)(i,j) = (\pi(i),\sigma(j)),
\qquad (\pi,\sigma) \in S_d \times S_d,
\qquad
i,j \in [d].
\]
Then the following hold:
\begin{itemize}
\item[(a)]
if some $(\pi,\sigma) \in S_d \times S_d$ induces a weak equivalence from $H$ to $H'$, then $\iota(\pi,\sigma)$ induces a strong equivalence from $\{x_{ij}'\}$ to $\{x_{ij}\}$;
\item[(b)]
if some $\pi \in S_{[d] \times [d]}$ induces a weak equivalence from $\{x_{ij}\}$ to $\{x_{ij}'\}$, then $\pi$ induces a strong equivalence from $\tilde{H}$ to $\tilde{H}'$;
\item[(c)]
$\iota(\operatorname{Aut}_w(H))
\leq \operatorname{Aut}_s(\{x_{ij}\}) 
\leq \operatorname{Aut}_w(\{x_{ij}\})
\leq \operatorname{Aut}_s(\tilde{H})$.
\end{itemize}
\end{theorem}

Part of our argument for Theorem~\ref{thm:sandwich} will be useful again later, so we isolate it as a lemma.

\begin{lemma}
\label{lem:omegas}
Abbreviating $X:=[d] \times [d]$, let $\{ x_i \}_{i\in X}$ and $\{ x_i' \}_{i \in X}$ be SICs constructed by Theorem~\ref{thm.main}.
If some $\pi \in S_X$ induces a weak equivalence from $\{x_i\}$ to $\{x_i'\}$, then there exist
    \begin{itemize}
        \item[(i)]
        $\varepsilon \in \{\pm1\}$ with $\varepsilon = 1$ if $\operatorname{char} K \neq 3$, 
        \item[(ii)]
        $\gamma \in \{ \text{id}, \overline{\cdot}\}$ with $\gamma = \text{id}$ if $\pi$ induces a strong equivalence from $\{x_i\}$ to $\{x_i'\}$, and
        \item[(iii)]
        $\omega_i \in \{\pm1,\pm\mathrm{i}\}$
    \end{itemize}
such that $\big(x_{\pi(i)}',x_{\pi(j)}'\big) = \varepsilon \omega_i \overline{\omega_j} (x_i,x_j)^\gamma$ for every $i,j \in X$.
\end{lemma}

\begin{proof}
We may assume $|X|>1$.
Suppose $\pi \in S_X$ induces a weak equivalence from $\{x_i\}$ to $\{x_i'\}$, and select $\gamma$, $U$, and $c_i$ as in Definition~\ref{defn:SIC equivalence}, so that $x_{\pi(i)}' = c_i (U x_i)^\gamma$, with $\gamma = \text{id}$ if $\pi$ induces a strong equivalence from $\{x_i\}$ to $\{x_i'\}$.
Define $\varepsilon$ to be the constant value of $\overline{c_i} c_i$.
To begin, we show that $\varepsilon \in \{\pm 1\}$, with $\varepsilon=1$ if $\operatorname{char} K \neq 3$.
For any $i \neq j$, we have
\begin{equation}
\label{eq:equivGram}
\big( x_{\pi(i)}', x_{\pi(j)}' \big) 
= \big( c_i (Ux_i)^\gamma, c_j (Ux_j)^\gamma \big) 
= \overline{c_i} c_j (U x_i, U x_j)^\gamma 
= \overline{c_i} c_j (x_i,x_j)^\gamma.
\end{equation}
Here, recall that both $\big( x_{\pi(i)}', x_{\pi(j)}' \big)$ and $(x_i,x_j)$ belong to $\{\pm 4, \pm 4 \mathrm{i}\}$ by Lemma~\ref{lem.gram matrix}.
Multiplying each side of~\eqref{eq:equivGram} by its conjugate, we find
\[
16 = \big( x_{\pi(i)}', x_{\pi(j)}' \big) \overline{ \big( x_{\pi(i)}', x_{\pi(j)}' \big) } 
= \overline{c_i} c_j (x_i,x_j)^\gamma \cdot c_i \overline{c_j} \overline{(x_i,x_j)^\gamma}
= 16 \overline{c_i} c_i \overline{c_j} c_j
= 16 \varepsilon^2.
\]
Since $\operatorname{char} K \neq 2$, it follows that $\varepsilon \in \{1,-1\}$.
Furthermore, just as in~\eqref{eq:equivGram},
\[
12 = \big( x_{\pi(i)}',x_{\pi(i)}' \big) 
= \overline{c_i} c_i (x_i,x_i)^\gamma
= 12 \varepsilon,
\]
so that $\varepsilon = 1$ if $\operatorname{char} K \neq 3$.

Next, fix some $k \in X$, and define $\omega_i := \overline{c_i} c_k$ for each $i \in X$; in particular, $\omega_k = c_k \overline{c_k} = \varepsilon \in \{\pm 1\}$.
More generally, for any $i \neq k$, substituting $j = k$ in~\eqref{eq:equivGram} and rearranging yields
\[
\omega_i = \overline{c_i} c_k 
= \frac{ \big( x_{\pi(i)}', x_{\pi(k)}' \big) }{ (x_i,x_k)^\gamma } 
\in \{ \pm 1, \pm \mathrm{i} \},
\]
where we have again used the fact that both $\big( x_{\pi(i)}', x_{\pi(k)}' \big)$ and $(x_i,x_k)^\gamma$ belong to $\{ \pm 4, \pm 4\mathrm{i} \}$.
Given any $i,j \in X$ we follow~\eqref{eq:equivGram} to obtain
\[
\big( x_{\pi(i)}', x_{\pi(j)}' \big) 
= \overline{c_i} c_j (x_i,x_j)^\gamma
= c_k^{-1} \omega_i \overline{c_k}^{-1} \overline{\omega_j} (x_i,x_j)^\gamma
= \varepsilon \omega_i \overline{\omega_j} (x_i,x_j)^\gamma.
\]
Finally, since $(x_i,x_j) \in \{ \pm 4, \pm 4 \mathrm{i}\}$ whenever $i \neq j$, we either have $(x_i,x_j)^\gamma = (x_i,x_j)$ for every $i,j$ (in the case where $\mathrm{i}^\gamma = \mathrm{i}$) or $(x_i,x_j)^\gamma = \overline{(x_i,x_j)}$ for every $i,j$ (in the case where $\mathrm{i}^\gamma = - \mathrm{i}$).
This completes the proof.
\end{proof}

We now apply Lemma~\ref{lem:omegas} to prove Theorem~\ref{thm:sandwich}.

\begin{proof}[Proof of Theorem~\ref{thm:sandwich}]
It suffices to prove (a) and (b).

For (a), suppose $(\pi,\sigma) \in S_d \times S_d$ induces a weak equivalence from $H$ to $H'$, and select the requisite $\varepsilon_i, \varepsilon'_i \in \{ \pm 1 \}$ such that $H_{\pi(i),\sigma(j)}' = \varepsilon_i \varepsilon'_j H_{ij}$ for every $i,j \in [d]$.
Consider the permutation matrices $P,Q\in K^{d\times d}$ given by $Pe_j=e_{\pi(j)}$ and $Qe_j=e_{\sigma(j)}$, and define $D := \operatorname{diag}(\varepsilon_1,\ldots,\varepsilon_d)$ and $D' := \operatorname{diag}(\varepsilon'_1,\ldots,\varepsilon'_d)$; thus, $P^\top H Q = D H' D'$.
We claim that
\[
x_{\pi(i),\sigma(j)}
=\varepsilon'_j PDx_{ij}',
\]
so that $\iota(\pi,\sigma)$ induces a strong equivalence from $\{x_{ij}'\}$ to $\{x_{ij}\}$.
Toward that end, let $h_j \in K^d$ denote the $j$th column of $H$, and let $h_j' \in K^d$ denote the $j$th colum of $H'$.
First,
\[
h_{\sigma(j)}
=He_{\sigma(j)}
=HQe_j
=PDH'D'e_j
=\varepsilon'_j PDH'e_j
=\varepsilon'_j PDh_j',
\]
and so
\[
x_{\pi(i),\sigma(j)}
=h_{\sigma(j)}\circ(\mathbf{1}+ze_{\pi(i)})
=(\varepsilon'_j PDh_j')\circ(\mathbf{1}+ze_{\pi(i)})
=\varepsilon'_j \big(PDh_j'+z(PDh_j')_{\pi(i)}e_{\pi(i)}\big).
\]
Next,
\[
(PDh_j')_{\pi(i)}e_{\pi(i)}
=(Dh_j')_iPe_i
=D_{ii}(h_j')_iPe_i
=(h_j')_iPDe_i
=PD(h_j'\circ e_i),
\]
and so
\begin{align*}
x_{\pi(i),\sigma(j)}
&=\varepsilon'_j \big(PDh_j'+z(PDh_j')_{\pi(i)}e_{\pi(i)}\big)\\
&=\varepsilon'_j \big(PDh_j'+zPD(h_j'\circ e_i)\big)
=\varepsilon'_j PD\big(h_j'\circ(\mathbf{1}+ze_i)\big)
=\varepsilon'_j PDx_{ij}',
\end{align*}
as claimed.

For (b), let $G$ and $G'$ denote the entrywise squares of the Gram matrices of $\{x_{ij}\}$ and $\{x_{ij}'\}$, respectively: $G_{(i,j),(k,l)} = (x_{ij},x_{kl})^2$ and $G_{(i,j),(k,l)}' = (x_{ij}',x_{kl}')^2$.
We will prove (b) by establishing two claims:
\begin{itemize}
\item[(i)]
if $\pi \in S_{[d] \times [d]}$ induces a weak equivalence from $\{x_{ij}\}$ to $\{x_{ij}'\}$, then it induces a strong equivalence from $G$ to $G'$;
\item[(ii)] 
if $\pi \in S_{[d] \times [d]}$ induces a strong equivalence from $G$ to $G'$, then it induces a strong equivalence from $\tilde{H}$ to $\tilde{H}'$.
\end{itemize}
For (i), suppose $\pi\in S_{[d] \times [d]}$ induces a weak equivalence from $\{ x_{ij} \}$ to $\{ x_{ij}' \}$.
Select $\gamma \in \{ \text{id}, \overline{\cdot} \}$, $\varepsilon \in \{\pm 1\}$, and $\omega_{ij} \in \{ \pm 1, \pm \mathrm{i} \}$ as in the conclusion of Lemma~\ref{lem:omegas}, so that 
\[
\big( x_{\pi(i,j)}', x_{\pi(k,l)}' \big) = \varepsilon \omega_{ij} \overline{ \omega_{kl} } (x_{ij}, x_{kl} )
\]
for every $(i,j), (k,l) \in [d] \times [d]$.
We square both sides and observe that $\overline{\omega_{kl}}^2 = \omega_{kl}^2$ and $\big( (x_{ij}, x_{kl})^\gamma \big)^2 = (x_{ij}, x_{kl})^2$ to obtain
\[
\big( x_{\pi(i,j))}', x_{\pi(k,l)}' \big)^2
=\omega_{ij}^2\omega_{kl}^2(x_{ij},x_{kl})^2.
\]
Since $\omega_{ij}^2 \in \{\pm1\}$, $\pi$ induces a strong equivalence from $G$ to $G'$.

For (ii), Lemma~\ref{lem.gram matrix} immediately implies that
\[
G_{(i,j),(k,l)}
=(x_{ij},x_{kl})^2
=16H_{ij}H_{il}H_{kj}H_{kl}
\qquad
\text{whenever}
\qquad
(i,j)\neq(k,l).
\]
Thus, all of the entries of $F:=\frac{1}{16}G-8I$ are given by
\[
F_{(i,j),(k,l)}
= H_{ij}H_{il}H_{kj}H_{kl}
= H_{ij}([r_k r_i^\top]_{jl})H_{kl}
= H_{ij} \tilde{H}_{(i,j), (k,l)} H_{kl}.
\]
Similarly, $F':=\frac{1}{16}G'-8I$ has entries
\[
F_{(i,j),(k,l)}' = H_{ij}' \tilde{H}_{(i,j), (k,l)}' H_{kl}'.
\]
It follows easily that any $\pi \in S_{[d] \times [d]}$ inducing a strong equivalence from $G$ to $G'$ also induces a strong equivalence from $\tilde{H}$ to $\tilde{H}'$.
\end{proof}

To exactly compute the symmetries of a SIC constructed by Theorem~\ref{thm.main}, we rely on the following proposition, the main idea of which is well known; see, for instance,~\cite{Cameron:77,OCathain:08}.

\begin{proposition}
\label{prop:AutCompute}
Abbreviating $X:=[d] \times [d]$, let $\{ x_i \}_{i \in X}$ be a SIC constructed by Theorem~\ref{thm.main}.
For each $\omega \in C_4:=\{\pm1,\pm \mathrm{i}\}$, consider the directed graph $\Gamma_\omega$ on vertex set $C_4 \times X$ with
\[
(\alpha,i) \to (\beta,j)
\text{ in }
\Gamma_\omega
\iff
(\alpha x_i, \beta x_j) = 4\omega.
\]
Then the following hold for any $\pi \in S_X$:   
    \begin{itemize}
    \item[(a)]
    $\pi \in \operatorname{Aut}_w(\{x_i\})$ if and only if there exist
        \begin{itemize}
        \item[(i)]
        a bijection $f\colon C_4 \times X \to C_4 \times X$ with $f(C_4 \times \{i\}) = C_4 \times \{ \pi(i)\}$ for every $i \in X$,
        \item[(ii)]
        $\varepsilon \in \{\pm1\}$ with $\varepsilon = 1$ if $\operatorname{char} K \neq 3$, and
        \item[(iii)]
        $\gamma \in \{\text{id},\overline{\cdot}\}$
        \end{itemize}
    such that $f$ induces an isomorphism $\Gamma_\omega \to \Gamma_{\varepsilon \omega^\gamma}$ for every $\omega \in C_4$;

    \item[(b)]
    $\pi \in \operatorname{Aut}_s(\{x_i\})$ if and only if there exist
        \begin{itemize}
        \item[(i)]
        a bijection $f\colon C_4 \times X \to C_4 \times X$ with $f(C_4 \times \{i\}) = C_4 \times \{ \pi(i)\}$ for every $i \in X$, and
        \item[(ii)]
        $\varepsilon \in \{\pm1\}$ with $\varepsilon = 1$ if $\operatorname{char} K \neq 3$
        \end{itemize}
    such that $f$ induces an isomorphism $\Gamma_\omega \to \Gamma_{\varepsilon \omega}$ for every $\omega \in C_4$.
    \end{itemize}

\end{proposition}

\begin{proof}
We prove (a), and (b) is similar.
In the forward direction, suppose $\pi \in \operatorname{Aut}_w(\{x_i\})$.
Apply Lemma~\ref{lem:omegas} with $x_i' = x_i$, and let $\varepsilon \in \{ \pm 1 \}$, $\gamma \in \{ \text{id}, \overline{\cdot} \}$, and $\omega_i \in C_4$ be as in the conclusion of the lemma, so that
\begin{equation}
\label{eq:weakOmegas}
(x_{\pi(i)}, x_{\pi(j)}) = \varepsilon \omega_i \overline{\omega_j} (x_i, x_j)^\gamma
\end{equation}
for every $i,j \in X$.
Given any $\alpha,\beta \in C_4$, we rearrange~\eqref{eq:weakOmegas} to obtain
\[
\big(\alpha^\gamma \omega_i x_{\pi(i)}, \, \beta^\gamma \omega_j x_{\pi(j)} \big) = \varepsilon\big(\alpha x_i, \beta x_j\big)^\gamma.
\]
Thus, for each $\omega \in C_4$,
\[
(\alpha,i) \to (\beta,j) \text{ in } \Gamma_\omega
\iff
\big(\alpha^\gamma \omega_i, \pi(i) \big) \to \big( \beta^\gamma \omega_j, \pi(j) \big) \text{ in } \Gamma_{\varepsilon \omega^\gamma}.
\]
The conclusion of (a) is therefore satisfied with $f \colon C_4 \times X \to C_4 \times X$ given by $f(\alpha,i) = \big( \alpha^\gamma \omega_i, \pi(i) \big)$.

Conversely, suppose there exist $f \colon C_4 \times X \to C_4 \times X$, $\varepsilon \in \{ \pm 1 \}$, and $\gamma \in \{ \text{id}, \overline{\cdot} \}$ as in the conclusion of (a).
We will produce a tuple $\{ c_i \}_{i \in X}$ of scalars for which $\overline{c_i} c_i$ is constant, as well as an invertible linear map $U \colon L^d \to L^d$ that preserves the sesquilinear form $(\cdot, \cdot)$, such that $x_{\pi(i)} = c_i (U x_i)^\gamma$ for every $i \in X$. 
For each $i \in X$, we apply (i) to define $\omega_i \in C_4$ by $f(1,i) = (\omega_i, \pi(i))$.
Given $i \neq j$ in $X$, Lemma~\ref{lem.gram matrix} provides $\omega \in C_4$ for which $(x_i, x_j) = 4\omega$, so that $(1,i) \to (1,j)$ in $\Gamma_\omega$.
Since $f$ induces an isomorphism $\Gamma_\omega \to \Gamma_{\varepsilon \omega^\gamma}$, it also holds that $(\omega_i, \pi(i) ) \to (\omega_j, \pi(j))$ in $\Gamma_{\varepsilon \omega^\gamma}$; that is,
\[
\big( \omega_i x_{\pi(i)}, \omega_j x_{\pi(j)} \big) 
= 4\varepsilon \omega^\gamma 
= \varepsilon(x_i,x_j)^\gamma.
\]
Rearranging, we find
\[
\big( x_{\pi(i)}, x_{\pi(j)} \big)
= \varepsilon \omega_i \overline{\omega_j} (x_i,x_j)^\gamma
\]
whenever $i \neq j$.
Furthermore, this equation holds even when $j = i$, as seen in cases: if $\varepsilon = 1$, then both sides equal~12; if $\varepsilon = -1$, then we have $\operatorname{char} K = 3$ by (ii), so both sides equal~0.
With this in mind, we define $\{ c_i \}_{i\in X}$ as follows: if $\varepsilon = 1$, then $c_i := \overline{\omega_i}$; if $\varepsilon = -1$, then $c_i := (1+\mathrm{i}) \overline{\omega_i}$.
Thus, $\overline{c_i} c_j = \varepsilon \omega_i \overline{\omega_j}$ and
\begin{equation}
\label{eq:unitary relation}
\big( x_{\pi(i)}, x_{\pi(j)} \big)
= \overline{c_i} c_j (x_i,x_j)^\gamma
\qquad
\text{for every }
i,j \in X.
\end{equation}
We will now apply~\eqref{eq:unitary relation} to produce the appropriate transformation $U \colon L^d \to L^d$ with $x_{\pi(i)} = c_i (U x_i)^\gamma$ for every $i \in X$; that is, $Ux_i = \big( c_i^{-1} x_{\pi(i)} \big)^\gamma =: y_i$.
For this purpose, it is convenient to rewrite~\eqref{eq:unitary relation} again as
\begin{equation}
\label{eq:unitary relation 2}
(y_i,y_j) = (x_i,x_j)
\qquad
\text{for every }
i,j \in X.
\end{equation}
To see the mapping $U x_i = y_i$ extends to well-defined linear transformation $U \colon L^d \to L^d$, first recall that $\operatorname{span}\{x_i : i \in X\} = L^d$.
Considering that spanning sets are precisely those containing bases, it follows easily that $\operatorname{span}\{y_i : i \in X \} = L^d$ as well.
Now suppose $\sum_i \alpha_i x_i = \sum_i \beta_i x_i$ for some scalars $\alpha_i,\beta_i \in L$.
Given any $j \in X$, we apply~\eqref{eq:unitary relation 2} to find
\[
\big( \sum_i (\alpha_i-\beta_i) y_i, y_j \big)
= \big( \sum_i (\alpha_i - \beta_i) x_i, x_j \big)
= 0.
\]
Since each $e_k \in L^d$ can be expressed as a linear combination of $\{y_j\}_{j\in X}$, it follows easily that $\sum_i (\alpha_i-\beta_i) y_i = 0$; that is, $\sum_i \alpha_i y_i = \sum_i \beta_i y_i$.
Overall, the function $U \colon L^d \to L^d$ with $U \sum_i \alpha_i x_i = \sum_i \alpha_i y_i$ is well defined.
It is clearly linear, and it is invertible since it maps a spanning set to a spanning set.
Finally, $U$ preserves the sesquilinear form by virtue of~\eqref{eq:unitary relation 2} and the fact that $\{x_i\}$ and $\{y_i\}$ both span $L^d$.
This completes the proof of (a).
The proof of (b) is nearly identical, the only modification being that $\gamma = \text{id}$ throughout.
\end{proof}

\section{Examples and discussion}

We now compute automorphism groups of SICs constructed from small Hadamard matrices in Theorem~\ref{thm.main}.
The order $d$ of the Hadamard matrix must satisfy $d \equiv 8$ in $K$, and so (with finite fields in mind) we restrict our attention to those $d$ for which $d-8$ is divisible by a prime $p \equiv 3 \bmod 4$.
The smallest such $d$ for which an honest Hadamard matrix exists are $d=2,8,20,32,36$.
For each of these values, we describe the symmetries of at least one SIC below.
In each case, we used the GRAPE package~\cite{GRAPE} of GAP~\cite{GAP} to find the requisite graph isomorphisms from Proposition~\ref{prop:AutCompute}.
Meanwhile, we computed automorphism groups of Hadamard matrices using standard techniques (such as in~\cite{OCathain:08}) and the DESIGN~\cite{DESIGN} and GRAPE packages of GAP.

\begin{example}
The Hadamard matrix $\left[ \begin{smallmatrix} 1 & 1 \\ 1 & -1 \end{smallmatrix} \right]$ of order~2 gives a SIC in~$\mathbb{F}_9^2$:
\[
\left[ \begin{array}{cccc}
x_{1,1} & x_{1,2} & x_{2,1} & x_{2,2}
\end{array} \right]
=
\left[ \begin{array}{cccc}
-1 - 2\mathrm{i} & -1 - 2\mathrm{i} & 1 & 1 \\
1 & -1 & -1 - 2\mathrm{i} & 1 + 2\mathrm{i}
\end{array} \right].
\]
This SIC also appeared in Theorem~6.5 of~\cite{GreavesIJM:22}, as part of an infinite family of finite field SICs with Weyl--Heisenberg symmetry.
Here, we find $\operatorname{Aut}_s(\{x_{ij}\}) = \operatorname{Aut}_w(\{x_{ij}\}) = S_{[2]\times [2]}$.
Half of the strong automorphisms (namely, the even permutations) arise by taking $\varepsilon = 1$ in Proposition~\ref{prop:AutCompute}(b), while the other half arise only by taking $\varepsilon = -1$.
Notably, $\{x_{ij}\}$ exhibits more symmetry than any complex SIC, since the strong automorphism group of a complex SIC cannot be 3-transitive~\cite{Zhu:15,King:19,IversonM:24}.
\end{example}

\begin{example}
\label{ex:Hoggar syms}
Returning to Example~\ref{ex:Hoggar}, consider Hoggar's lines in $\mathbb{C}^8$, which arise from a Hadamard matrix of order~8.
Hoggar's lines are highly symmetric, and they are known to give one of exactly three complex SICs whose strong automorphism group is 2-transitive~\cite{Zhu:15,DempwolffK:23,IversonM:24}.
In the automorphism group sandwich of Theorem~\ref{thm:sandwich}(c), we find all three inclusions are strict.
Specifically:
\begin{itemize}
    \item 
    $\iota(\operatorname{Aut}_w(H))\cong \mathbb{F}_2^3 \rtimes (\mathbb{F}_2^3 \rtimes \operatorname{SL}(3,2))$ is transitive but not 2-transitive,
    \item
    $\operatorname{Aut}_s(\{x_{ij}\}) \cong \mathbb{F}_2^6 \rtimes G_2(2)'$ is 2-transitive and contains $\iota(\operatorname{Aut}_w(H))$ with index~36,
    \item
    $\operatorname{Aut}_w(\{x_{ij}\})$ contains $\operatorname{Aut}_s(\{x_{ij}\})$ with index~2, and
    \item
    $\operatorname{Aut}_s(\tilde{H}) \cong \mathbb{F}_2^6 \rtimes \operatorname{Sp}(6,2)$ contains $\operatorname{Aut}_w(\{x_{ij}\})$ with index~120.
\end{itemize}
\end{example}

\begin{example}
Up to weak equivalence, there are exactly three Hadamard matrices of order~20~\cite{Hall:65}.
Each gives a SIC in $\mathbb{F}_9^{20}$.
In all three cases, the weak automorphism group of the SIC is not transitive.
For the ``Paley~I'' construction, the orbits have sizes 380 and~20;
for the ``Paley~II'' construction, the orbit sizes are 360 and~40;
for the remaining construction, the orbit sizes are 320 and~80.
For all three constructions, the automorphism group sandwich of Theorem~\ref{thm:sandwich}(c) has exactly two distinct layers, where $\iota(\operatorname{Aut}_w(H))=\operatorname{Aut}_s(\{x_{ij}\})$ has index~2 in $\operatorname{Aut}_w(\{x_{ij}\})=\operatorname{Aut}_s(\tilde{H})$.
\end{example}

\begin{example}
For each prime power $q \equiv 3 \bmod 4$, the ``Paley I'' construction creates a Hadamard matrix $H$ of order $d:=q+1$, and for $q>11$ it is known that $\operatorname{Aut}_w(H) \cong \operatorname{P\Sigma L}(2,q)$~\cite{Kantor:69}.
If it happens that $d \equiv 8 \bmod p$ for some prime $p \equiv 3 \bmod 4$, then Theorem~\ref{thm.main} creates a SIC $\{x_{ij}\}$ in $\mathbb{F}_{p^2}^d$.
For each of the computationally feasible values $q = 19,31,43,67$, we observe that $\operatorname{Aut}_s(\{x_{ij}\}) = \operatorname{Aut}_w(H) \cong \operatorname{P\Sigma L}(2,q)$ and $\operatorname{Aut}_w(\{x_{ij}\}) = \operatorname{Aut}_s(\tilde{H}) \cong \operatorname{P\Gamma L}(2,q)$.
In particular, $\operatorname{Aut}_w(\{x_{ij}\})$ is not transitive, instead having orbits $\{(i,i) : i \in [d]\}$ and its complement.
Furthermore, for each of $q=23,27,47,59,71$ (where no appropriate prime $p$ exists) we continue to observe that $\operatorname{Aut}_s(\tilde{H}) \cong \operatorname{P\Gamma L}(2,q)$.
\end{example}

We expect this behavior to generalize.

\begin{conjecture}
\label{conj:Paley}
Given a prime power $q > 11$ with $q \equiv 3 \bmod 4$, consider the ``Paley I'' Hadamard matrix $H$ of order $d:=q+1$.
Then the following hold:
    \begin{itemize}
        \item[(a)]
        $\operatorname{Aut}_s(\tilde{H}) \cong P\Gamma L(2,q)$ has exactly two orbits on $[d] \times [d]$, namely, $\{(i,i) : i \in [d]\}$ and its complement, and
        \item[(b)]
        for any prime $p \equiv 3 \bmod 4$ that divides $d-8$, the SIC $\{x_{ij}\}$ in $\mathbb{F}_{p^2}^d$ constructed by Theorem~\ref{thm.main} has $\operatorname{Aut}_s(\{x_{ij}\}) = \iota(\operatorname{Aut}_w(H))$ and $\operatorname{Aut}_w(\{x_{ij}\}) = \operatorname{Aut}_s(\tilde{H})$. 
    \end{itemize}
In particular, there exist non-transitive SICs in infinitely many dimensions over various finite fields. 
\end{conjecture}

As progress toward Conjecture~\ref{conj:Paley}(a), the authors are able to show that $\operatorname{Aut}_s(\tilde{H})$ contains a copy of $P\Gamma L(2,q)$, whose orbits on $[d] \times [d]$ are exactly as described above.
(We omit our proof for the sake of brevity.)
However, we have not been able to prove the reverse inclusion of groups.
Part of the problem is that we do not see any transitive action of $\operatorname{Aut}_s(\tilde{H})$ to leverage as was done in~\cite{Kantor:69,deLauneyS:08}.
In particular, we do not know how to prove that $\operatorname{Aut}_s(\tilde{H})$ holds the diagonal $\{(i,i) : i \in [d]\}$ invariant.

\begin{example}
The homepage of E.\ Spence contains a file with 24 ``Goethals--Seidel''-type Hadamard matrices of order~36 that were catalogued by Spence and Turyn; of these 24 matrices, three have trivial weak automorphism group~\cite{SpenceT,Spence}.
Each of these three matrices creates a SIC in $\mathbb{F}_{49}^{36}$, and considering the automorphism group sandwich in Theorem~\ref{thm:sandwich}(c), one might expect these SICs to have few symmetries.
In fact, the SIC arising from matrix \#3 in~\cite{SpenceT} has trivial strong automorphism group, and its weak automorphism group has order~2.
Furthermore, each of matrix \#23 and matrix \#24 produces a \textbf{totally asymmetric} SIC, meaning its weak automorphism group is trivial.
\end{example}

In light of the automorphism group sandwich in Theorem~\ref{thm:sandwich}(c), and given the empirical abundance of asymmetric Hadamard matrices, we expect totally asymmetric SICs to exist commonly in the finite field setting, as follows.

\begin{conjecture}
There exist totally asymmetric SICs in infinitely many dimensions over various finite fields.
Specifically, for infinitely many choices of $d$, there exists a modular Hadamard matrix $H$ of order~$d$ for which $\operatorname{Aut}_s(\tilde{H})$ is trivial, and $d \equiv 8 \bmod p$ for some prime $p \equiv 3 \bmod 4$, so that the SIC in $\mathbb{F}_{p^2}^d$ constructed by Theorem~\ref{thm.main} has trivial weak automorphism group.
\end{conjecture}

We end with a brief discussion.
To date, the bulk of research on Zauner's conjecture has focused on complex SICs with Weyl--Heisenberg symmetry, and for good reason: the data suggests they exist in every dimension.
However, given the stubborn number-theoretic difficulties behind this approach, and given the existence of asymmetric SICs over finite fields, it might also be worthwhile searching for complex SICs whose existence relies more on combinatorics than on symmetry.
For example, Example~5.11 of~\cite{IversonM:22} identifies an infinite family of feasible parameters of \textit{roux} that would produce complex SICs from certain unknown graphs.
The first member of this family exists and produces Hoggar's lines in $\mathbb{C}^8$; the second (unknown) member would produce a SIC in $\mathbb{C}^{24}$, but the roux would necessarily lack certain symmetries.
(Specifically, the corresponding \textit{roux scheme} could not arise from a multiplicity-free permutation group.)
More generally, we pose the following problem, in three levels of intensity.

\begin{problem}
Does there exist a complex SIC\ldots
    \begin{itemize}
    \item[(a)]
    \ldots without Weyl--Heisenberg symmetry?
    \item[(b)]
    \ldots whose strong automorphism group is not transitive?
    \item[(c)]
    \ldots whose weak automorphism group is trivial?
    \end{itemize}
\end{problem}

We note that in dimension~$2$, there is a unique complex SIC up to strong equivalence, which in turn exhibits Weyl--Heisenberg symmetry.
Next, it is known~\cite{HughstonS:16} that every complex SIC in dimension~$3$ is strongly equivalent to one with Weyl--Heisenberg symmetry.
In dimension~$8$, Hoggar's lines have a different type of Weyl--Heisenberg symmetry than those conjectured to exist by Zauner, since the underlying abelian group $(\mathbb{Z}/2\mathbb{Z})^3$ is not cyclic.
Meanwhile, in prime dimensions, it is known~\cite{Zhu:10} that every SIC with a transitive strong automorphism group is strongly equivalent to a SIC with Weyl--Heisenberg symmetry.

\section*{Acknowledgments}

JWI was supported by NSF DMS
2220301 and a grant from the Simons Foundation.  
DGM was supported by NSF DMS 2220304.

\end{document}